# A Comparative Study of State Transition Algorithm with Harmony Search and Artificial Bee Colony

Xiaojun Zhou[1,2], David Yang Gao[1], Chunhua Yang[2]

[1] School of Science, Information Technology and Engineering, University of Ballarat,
Victoria 3350, Australia
[2] School of Information Science and Engineering, Central South University,
Changsha 410083, China
tiezhongyu2010@gmail.com

**Abstract.** We focus on a comparative study of three recently developed nature-inspired optimization algorithms, including state transition algorithm, harmony search and artificial bee colony. Their core mechanisms are introduced and their similarities and differences are described. Then, a suit of 27 well-known benchmark problems are used to investigate the performance of these algorithms and finally we discuss their general applicability with respect to the structure of optimization problems.

## 1   Introduction

Existing natural phenomena, such as natural selection and survival of the fittest (genetic algorithm), natural annealing process in metallurgy (simulated annealing), foraging behavior of real ant colonies (ant colony optimization), and social behavior of bird flocks and fish schools (particle swarm optimization) have inspired researchers to develop algorithms for optimization problems. These nature-inspired algorithms have received considerable attention due to their strong adaptability and easy implementation. Inspired by the improvisation process of musicians and foraging behavior of real honeybees, harmony search (HS) [1-3] and artificial bee colony (ABC) [4-5] have been proposed respectively in recent few years. At the same time, in terms of the concepts of state and state transition, a new heuristic random search algorithm named state transition algorithm (STA) has been introduced in order to probe into classical and intelligent optimization algorithms [6-9]. In STA, a solution to a specific optimization problem is regarded as a state, and the process to update current state is considered as a state transition. For continuous optimization problem, it has designed four special geometrical operators called rotation, translation, expansion and axesion for both local exploitation and global exploration. For discrete optimization problem, some intelligent operators such as swap, shift, symmetry and substitute are designed for adapting to various types of problems. In this study, we focus on a comparative study of state transition algorithm with harmony search and artificial bee colony in their standard versions.

## 2 Three Stochastic Algorithms

In this section, we give a brief description of the three stochastic algorithms with respect to their mechanisms, and the similarities and differences are also discussed.

### 2.1 Harmony Search

In HS, there exist three possible choices to generate a new piece of music: (1) select a note stored in harmony memory at a probability of HMCR (harmony memory considerate rate); (2) adjust the pitch slightly at a probability of PAR (pitch adjusting rate); (3) compose any pitch randomly within bounds. The pitch is adjusted by

$$x_{new} = x_{old} + (2rand - 1) * b$$

where, $rand$ is a random number from [0,1], and $b$ is the pitch bandwidth.

### 2.2 Artificial Bee Colony

In ABC, the colony of artificial bees contains three groups of bees: (1) employed bees, going to the food source visited previously; (2) onlookers, making decision to choose a food source; (3) scouts, carrying out random search. A new position is produced by

$$v_{ij} = x_{ij} + \phi_{ij}(x_{ij} - x_{kj})$$

where, $i$ is the index of $i$th food position, $j$ is the $j$th component of a position, $k$ is a different index from $i$, and $j$, $k$ are created randomly.

An artificial onlooker bee chooses a food source depending on a probability by

$$p_i = \frac{\text{fit}_i}{\sum_{n=1}^{SN} \text{fit}_n}$$

where, $\text{fit}_i$ is the fitness value of the position $i$, $SN$ is the number of food sources.

### 2.3 State Transition Algorithm

In STA, there are four special geometrical operators defined by
(1) Rotation transformation

$$x_{k+1} = x_k + \alpha \frac{1}{n\|x_k\|_2} R_r x_k,$$

where, $\alpha$ is a positive constant, $R_r$ is a random matrix with its entries from [-1,1].
(2) Translation transformation

$$x_{k+1} = x_k + \beta R_t \frac{x_k - x_{k-1}}{\|x_k - x_{k-1}\|_2},$$

where, $\beta$ is a positive constant, $R_t$ is a random variable from [0,1].
(3) Expansion transformation

$$x_{k+1} = x_k + \gamma R_e x_k,$$

where, $\gamma$ is a positive constant, $R_e$ is a random diagonal matrix with its entries obeying the standard norm distribution.

(4) Axesion transformation

$$x_{k+1} = x_k + \delta R_a x_k$$

where, $\delta$ is a positive constant, $R_a$ is a random diagonal matrix with its entries obeying the standard norm distribution and only one random position having nonzero value.

### 2.4 Similarities and Difference

There are two main similarities among the three algorithms in the discussed versions: Firstly, a new solution is created randomly, and they are all stochastic algorithms. Second, "greedy criterion" is adopted to evaluate a solution, and it is different from simulated annealing, in which, a bad solution is accepted in probability.

The differences between STA and other two algorithms are: (1) both HS and ABC focus on updating each component of a solution, while STA treats a solution in whole for update except the axesion transformation; (2) the comparing STA is individual-based, while both HS and ABC are population-based; (3) the mutant operators are different in three algorithms; (4) in HS, there is a probability in choosing an update, while in STA, the updating procedures are determined; (5) in ABC, choosing a food source depending on a probability associated with the fitness, while in STA, a candidate solution with better fitness is preferred; (6) in ABC, the fitness is standardized, while in STA, the fitness is based on objective function.

## 3 Experimental Results

All these benchmark instances are taken from [10]. In our experiments, we use the codes of standard HS and ABC from [11, 12], and comparing version of the STA is from [7]. The size of the population is 10, and the maximum iterations (make sure that the maximum number of function evaluations is the same) are 1e3, 2e3, 4e3, 1e4, 5e4, and 1e5 for $n = 2, 3, 4, 10, (20, 24, 25)$ and 30, respectively. For each benchmark instance, the initial population is the same for three algorithms at a run, and 20 runs are performed for each algorithm. Statistics like mean, std (standard deviation), and Wilcoxon rank sum test are used to evaluate the STA with other two algorithms.

### 3.1 Benchmark instances

The details of the benchmark instances are given as follows.
Ackley function

$$f_1(x) = -20\exp(-0.2\sqrt{\frac{1}{n}\sum_{i=1}^{n}x_i^2}) - \exp(\frac{1}{n}\sum_{i=1}^{n}\cos(2\pi x_i)) + 20 + e, -15 \leq x_i \leq 30$$

Beale function
$$f_2(x) = (1.5 - x_1 + x_1 x_2)^2 + (2.25 - x_1 + x_1 x_2^2)^2 + (2.625 - x_1 + x_1 x_2^3)^2, -4.5 \leq x_i \leq 4.5$$

Bohachevsky Function
$$f_3(x) = x_1^2 + 2x_2^2 - 0.3\cos(3\pi x_1) - 0.4\cos(4\pi x_2) + 0.7, -100 \leq x_i \leq 100$$

Booth Function $f_4(x) = (x_1 + 2x_2 - 7)^2 + (2x_1 + x_2 - 5)^2, -10 \leq x_i \leq 10$

Branin Function
$$f_5(x) = (x_2 - \frac{5.1}{4\pi^2}x_1^2 + \frac{5}{\pi}x_1 - 6)^2 + 10(1 - \frac{1}{8\pi})\cos(x_1) + 10, -5 \leq x_1 \leq 10, 0 \leq x_2 \leq 15$$

Colville Function
$$f_6(x) = 100(x_1^2 - x_2)^2 + (x_1 - 1)^2 + (x_3 - 1)^2 + 90(x_3^2 - x_4)^2 + 10.1((x_2 - 1)^2 + (x_4 - 1)^2)$$
$$+ 19.8(x_2 - 1)(x_4 - 1), -10 \leq x_i \leq 10$$

Dixon & Price Function $f_7(x) = (x_1 - 1)^2 + \sum_{i=2}^{n} i(2x_i^2 - x_{i-1})^2, -10 \leq x_i \leq 10$

Easom Function
$$f_8(x) = -\cos(x_1)\cos(x_2)\exp(-(x_1 - \pi)^2 - (x_2 - \pi)^2), -100 \leq x_i \leq 100$$

Goldstein & Price Function
$$f_9(x) = (1 + (x_1 + x_2 + 1)^2 (19 - 14x_1 + 13x_1^2 - 14x_2 + 6x_1 x_2 + 3x_2^2)) \times$$
$$(30 + (2x_1 - 3x_2)^2 (18 - 32x_1 + 12x_1^2 - 48x_2 - 36x_1 x_2 + 27x_2^2)), -2 \leq x_i \leq 2$$

Griewank Function
$$f_{10}(x) = \frac{1}{4000}\sum_{i=1}^{n} x_i^2 - \prod_{i=1}^{n} \cos\left|\frac{x_i}{\sqrt{i}}\right| + 1, -600 \leq x_i \leq 600$$

Hartmann Function
$$f_{11}(x) = -\sum_{i=1}^{4} a_i \exp\left[-\sum_{j=1}^{3} A_{ij}(x_j - P_{ij})^2\right], 0 < x_j < 1$$

where,

$$a = [1, 1.2, 3, 3.2]^T, A = \begin{bmatrix} 3.0 & 10 & 30 \\ 0.1 & 10 & 35 \\ 3.0 & 10 & 30 \\ 0.1 & 10 & 35 \end{bmatrix}, P = 10^{-4}\begin{bmatrix} 6890 & 1170 & 2673 \\ 4699 & 4387 & 7470 \\ 1091 & 8732 & 5547 \\ 381 & 5743 & 8828 \end{bmatrix}$$

Hump Function $f_{12}(x) = 4x_1^2 - 2.1x_1^4 + \frac{1}{3}x_1^6 + x_1 x_2 - 4x_2^2 + 4x_2^4, -5 \leq x_i \leq 5$

Levy Function
$$f_{13}(x) = \sin^2(\pi y_1) + \sum_{i=1}^{n-1}\left[(y_i - 1)^2 (1 + 10\sin^2(\pi y_i + 1)) + (y_n - 1)^2 (1 + 10\sin^2(\pi y_n))\right]$$
$$y_i = 1 + \frac{x_i - 1}{4}, -10 \leq x_i \leq 10$$

Matyas Function $f_{14}(x) = 0.26(x_1^2 + x_2^2) - 0.48x_1x_2, -10 \leq x_i \leq 10$

Michalewics Function $f_{15}(x) = -\sum_{i=1}^{2} \sin(x_i) \sin(ix_i^2/\pi)^{2m}, m = 10, 0 \leq x_i \leq \pi$

Perm Functions $f_{16}(x) = \sum_{k=1}^{n}\left[\sum_{i=1}^{n}(i^k + \beta)\left((x_i/i)^k - 1\right)\right]^2, \beta = 0.5, -n \leq x_i \leq n$

Powell Function
$$f_{17}(x) = \sum_{i=1}^{n/4}(x_{4i-3} + 10x_{4i-2})^2 + 5(x_{4i-1} - x_{4i})^2 + (x_{4i-2} - x_{4i-1})^4 + 10(x_{4i-3} - x_{4i})^4, -4 \leq x_i \leq 5$$

Power Sum Function $f_{18}(x) = \sum_{k=1}^{n}\left[\left(\sum_{i=1}^{n}x_i^k\right) - b_k\right]^2, b = (8, 18, 44, 114), -4 \leq x_i \leq 5$

Rastrigin Function $f_{19}(x) = \sum_{i=1}^{n}(x_i^2 - 10\cos(2\pi x_i) + 10), -5.12 \leq x_i \leq 5.12$

Rosenbrock Function $f_{20}(x) = \sum_{i=1}^{n}(100(x_{i+1} - x_i^2)^2 + (x_i - 1)^2), -5 \leq x_i \leq 10$

Schwefel Function $f_{21}(x) = 418.9829n - \sum_{i=1}^{n}\left(x_i \sin\sqrt{|x_i|}\right), -500 \leq x_i \leq 500$

Shekel Function
$$f_{22}(x) = -\sum_{j=1}^{m}\left[\sum_{i=1}^{4}(x_i - C_{ij})^2 + \beta_j\right]^{-1}, m = 10, 0 \leq x_i \leq 10$$

$$\beta = \frac{1}{10}[1, 2, 2, 4, 4, 6, 3, 7, 5, 5]^T, C = \begin{bmatrix} 4.0 & 1.0 & 8.0 & 6.0 & 3.0 & 2.0 & 5.0 & 8.0 & 6.0 & 7.0 \\ 4.0 & 1.0 & 8.0 & 6.0 & 7.0 & 9.0 & 5.0 & 1.0 & 2.0 & 3.6 \\ 4.0 & 1.0 & 8.0 & 6.0 & 3.0 & 2.0 & 3.0 & 8.0 & 6.0 & 7.0 \\ 4.0 & 1.0 & 8.0 & 6.0 & 7.0 & 9.0 & 3.0 & 1.0 & 2.0 & 3.6 \end{bmatrix}$$

Shubert Function $f_{23}(x) = \sum_{i=1}^{5}i\cos[(i+1) \cdot x_1 + i] \cdot \sum_{i=1}^{5}i\cos[(i+1) \cdot x_2 + i], -10 \leq x_i \leq 10$

Sphere Function $f_{24}(x) = \sum_{i=1}^{n}x_i^2, -5.12 \leq x_i \leq 5.12$

Sum Squares Function $f_{25}(x) = \sum_{i=1}^{n}ix_i^2, -10 \leq x_i \leq 10$

Trid Function $f_{26}(x) = \sum_{i=1}^{n}(x_i - 1)^2 - \sum_{i=2}^{n}x_i x_{i-1}, -n^2 \leq x_i \leq n^2$

Zakharov Function
$$f_{27}(x) = \sum_{i=1}^{n}x_i^2 + \left(\sum_{i=1}^{n}0.5ix_i\right)^2 + \left(\sum_{i=1}^{n}0.5ix_i\right)^4, -5 \leq x_i \leq 10$$

## 3.2 Results and Discussion

Test results are listed in Table 1. We can find that the results of HS are always not as good as that of ABC and STA, except for $f_{11}$, $f_{15}$ and $f_{23}$. It seems that HS are capable of solving problems without much interaction between variables, and the solution accuracy and global search ability of HS are also not satisfactory. Considering that both HS and ABC focus on the mutation of each component of a solution, we can observe that, in HS, the pitch adjustment is a little blind, while in ABC, a new solution is created by adding a certain difference between current solution and a randomly different solution to current solution, which can be viewed as a strategy to share information and it is quite significant in swarm intelligence.

**Table 1.** Results for three algorithms on benchmark instances

| Functions | HS mean$\pm$std | ABC mean$\pm$std | STA mean$\pm$std |
|---|---|---|---|
| $f_1(n=2)$ | 0.14$\pm$0.57 $-$ | 8.88E-16$\pm$0 $\approx$ | 8.88E-16$\pm$0 |
| $f_2(n=2)$ | 0.35$\pm$0.53 $-$ | 3.61E-06$\pm$1.38E-5 $-$ | 4.31E-11$\pm$4.91E-11 |
| $f_3(n=2)$ | 0.73$\pm$0.62 $-$ | 0$\pm$0$\approx$ | 0$\pm$0 |
| $f_4(n=2)$ | 0.08$\pm$0.14 $-$ | 4.57E-17$\pm$4.90 E-17 + | 4.80E-11 $\pm$3.99E-11 |
| $f_5(n=2)$ | 0.39$\pm$0.01 $-$ | 0.39$\pm$5.46E-16 $\approx$ | 0.39$\pm$1.50E-16 |
| $f_6(n=4)$ | 7.20$\pm$19.98 $-$ | 0.21$\pm$0.14 $-$ | 0.001$\pm$0.002 |
| $f_7(n=25)$ | 12.17$\pm$5.19 $-$ | 7.51E-15$\pm$2.71E-15 + | 0.60$\pm$0.20 |
| $f_8(n=2)$ | -0.43$\pm$0.49 $-$ | -0.9057$\pm$0.27 $-$ | -1.0$\pm$1.31E-11 |
| $f_9(n=2)$ | 11.48$\pm$12.86 $-$ | 3.002$\pm$0.008 $-$ | 3.00$\pm$4.77E-9 |
| $f_{10}(n=2)$ | 0.16$\pm$0.14 $-$ | 0$\pm$0$\approx$ | 0$\pm$0 |
| $f_{11}(n=3)$ | -3.86$\pm$2.88E-8 $\approx$ | -3.86$\pm$1.82E-15$\approx$ | -3.86$\pm$2.96E-10 |
| $f_{12}(n=2)$ | 2.23E-5$\pm$8.82E-5$-$ | 4.65E-8$\pm$0 $\approx$ | 4.66E-8$\pm$1.13E-10 |
| $f_{13}(n=30)$ | 0.90$\pm$0.22 $-$ | 4.98E-16$\pm$5.39E-17 + | 3.84E-11$\pm$4.80E-12 |
| $f_{14}(n=2)$ | 0.05$\pm$0.06$-$ | 4.27E-10$\pm$1.75E-9 $-$ | 1.97E-250$\pm$0 |
| $f_{15}(n=2)$ | -1.8013$\pm$5.44E-5 $\approx$ | -1.8013$\pm$6.83E-16 $\approx$ | -1.8013 $\pm$1.01E-10 |
| $f_{16}(n=4)$ | 5.94$\pm$9.22 $-$ | 0.15$\pm$0.14 $-$ | 0.01$\pm$0.03 |
| $f_{17}(n=24)$ | 10.27$\pm$5.54 $-$ | 1.88E-4$\pm$5.94E-5$\approx$ | 1.13E-4$\pm$2.36E-5 |
| $f_{18}(n=4)$ | 0.29$\pm$0.49 $-$ | 0.02$\pm$0.01$-$ | 4.33E-4$\pm$5.02E-4 |
| $f_{19}(n=2)$ | 0.09$\pm$0.30 $-$ | 0$\pm$0$\approx$ | 0$\pm$0 |
| $f_{20}(n=2)$ | 1.02$\pm$1.49 $-$ | 0.01$\pm$0.01$-$ | 4.38E-8$\pm$1.71E-7 |
| $f_{21}(n=2)$ | 0.03$\pm$0.16 $-$ | 2.54E-5$\pm$0$\approx$ | 2.54E-5$\pm$1.48E-12 |
| $f_{22}(n=4)$ | -5.61$\pm$3.41 $-$ | -10.53$\pm$9.34E-5 $\approx$ | -10.53$\pm$3.06E-10 |
| $f_{23}(n=2)$ | -186.73$\pm$5.09E-4 $\approx$ | -186.73$\pm$3.57E-14 $\approx$ | -186.73$\pm$3.28E-8 |
| $f_{24}(n=30)$ | 0.72$\pm$0.20 $-$ | 5.08E-16$\pm$5.69E-17$-$ | 0$\pm$0 |
| $f_{25}(n=20)$ | 0.69$\pm$0.50 $-$ | 2.58E-16$\pm$3.72E-17$-$ | 0$\pm$0 |
| $f_{26}(n=10)$ | -78.37$\pm$113.69 $-$ | -210$\pm$7.32E-7$\approx$ | -210$\pm$1.86E-10 |
| $f_{27}(n=2)$ | 6.37E-4$\pm$2.80E-3$-$ | 2.91E-18$\pm$2.55E-18$-$ | 0$\pm$0 |

For ABC and STA, we can find their results are much more satisfactory, and they are able to obtain the global solutions for the majority of the test problems. To be more specific, we can find that ABC outperforms STA for $f_4$, $f_7$ and $f_{13}$, and it can gain higher precision than STA, especially for $f_7$, which indicates that ABC are more suitable for problems with strongly interacted structure. On the other hand, for $f_2$, $f_6$, $f_8$, $f_9$, $f_{14}$, $f_{16}$, $f_{18}$, $f_{20}$, $f_{24}$, $f_{25}$ and $f_{27}$, STA outperforms ABC in terms of solution accuracy, which indicates STA has stronger local exploitation ability than that of ABC.

Fig.1. gives the average fitness curve of Matyas function by the three algorithms. We can find that STA is more capable of searching in depth. The reason can be explained by the rotation transformation. In standard STA, a rotation factor is decreasing from a certain constant to an extremely small one in a periodical way.

However, due to the simple dimensional search along each axes and the base on individual search, standard STA are not good for problems with strongly interacted variables. Although population-based STA has been proposed in [8] and it indicated that communication can ameliorate its performance greatly, how to share information among different individuals are key issues. Furthermore, by this study, it shows that highlighted mutation in each component of a solution is a beneficial method, and these are useful indications for future development of state transition algorithm.

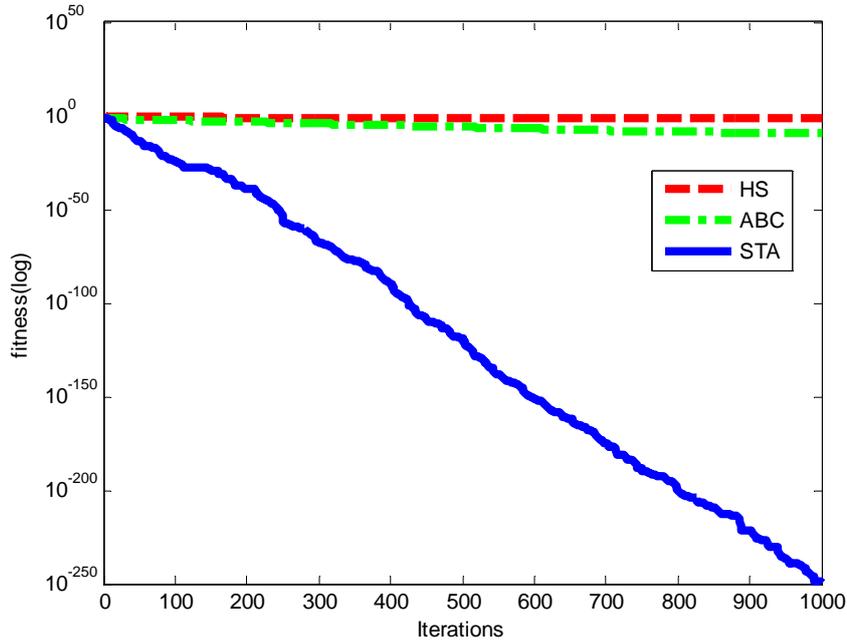

**Fig. 1. The average fitness curve of Matyas function by HS, ABC and STA**

## 4      Conclusion

In this paper, we investigate the mechanisms and performances of state transition algorithm, harmony search and artificial bee colony. Similarities and differences of the algorithms are mainly focused. A suit of unconstrained optimization problems have been used to evaluate these algorithms. Experimental results show that both state transition algorithm and artificial bee colony have better global search capability and can achieve higher solution accuracy than harmony search, artificial bee colony is more capable of solving problems with strongly interacted variables, and state transition algorithm has the potential ability to search in depth.

**Acknowledgments.** Xiaojun Zhou's research is supported by China Scholarship Council, David Yang Gao is supported by US Air Force Office of Scientific Research under the grant AFOSR FA9550-10-1-0487 and Chunhua Yang is supported by the National Science Found for Distinguished Young Scholars of China (Grant No.61025015).

## References

1. Geem, Z.W., Kim, J.H. and Loganathan, G.V.: A new heuristic optimization algorithm: Harmony Search. Simulation. 76(2)(2001) , 60-68
2. Lee, K.S., Geem, Z.W.: A new meta-heuristic algorithm for continuous engineering optimization: harmony search theory and practice. Comput. Methods Appl. Mech. Engrg. 194(2005) 3902-3933
3. Yang X.S. Harmony search as a metaheuristic algorithm. In: Music-Inspired Harmony Search Algorithm: Theory and Application, Springer, (2009) pp. 1-14
4. Karaboga, D., Basturk, B.: A powerful and efficient algorithm for numerical function optimization: artificial bee colony (ABC) algorithm. J. Glob. Optim. 39(2007) 459–471
5. Karaboga, D., Akay, B.: A comparative study of artificial bee colony algorithm. Appl. Math. Comput. 214(2009) 108–132
6. Zhou, X.J., Yang C.H., Gui, W.H.: Initial version of state transition algorithm. In Second International Conference on Digital Manufacturing and Automation (ICDMA), (2011), 644–647
7. Zhou, X.J., Yang C.H., Gui, W.H.: A new transformation into State Transition Algorithm for finding the global minimum. In 2nd International Conference on Intelligent Control and Information Processing (ICICIP), (2011), 674 –678
8. Zhou, X.J., Yang C.H., Gui, W.H.: State transition algorithm. Journal of Industrial and Management Optimization 8(4) (2012), 1039–1056
9. Yang, C.H., Tang, X.L., Zhou, X.J., Gui, W.H.: State transition algorithm for traveling salesman problem. to be published in the 31$^{st}$ Chinese Control Conference, arXiv:1206.0329
10.    http://www-optima.amp.i.kyotou.ac.jp/member/student/hedar/Hedar_files/TestGO_files/Page364
11.    https://sites.google.com/a/hydroteq.com/www/
12.    http://mf.erciyes.edu.tr/abc/